\documentclass[12pt]{article}
\usepackage{bbm}
\usepackage{latexsym}
\usepackage{mathrsfs}
\usepackage{graphicx}
\usepackage{subfigure}
\usepackage{amssymb}
\usepackage{amsmath}
\usepackage{amsthm}
\usepackage{color}
\usepackage{cite}
\usepackage{indentfirst}
\usepackage{anysize}\marginsize{45mm}{45mm}{40mm}{50mm}
\usepackage{caption}
\usepackage{float}
\usepackage{multirow}

\usepackage{enumerate}

\usepackage{latexsym, bm}
\newtheoremstyle{lemma}{\topsep}{\topsep}%
     {}
     {}
     {\bfseries}
     {}
     {0.1em}
     {\thmname{#1}\thmnumber{ #2}\thmnote{ #3}}
\theoremstyle{lemma}  

\newtheorem{theorem}{Theorem}     
\newtheorem{lemma}[theorem]{Lemma}
\newtheorem{corollary}[theorem]{Corollary}

\newtheorem{definition}{Definition}
\newtheorem{example}[theorem]{Example}

\numberwithin{equation}{section}

\title{Symmetric properties and two variants of shuffle-cubes\thanks{This research was partially supported by NSFC (Nos. 11801061 and 12161002). }}

\author{ Huazhong L\"{u}$^{1}$\thanks{Corresponding author.}, Kai Deng$^{2}$ and Xiaomei Yang$^{3}$\\
{\small $^{1}$School of Mathematical Sciences,} \\
{\small University of Electronic Science and Technology of China,}\\
{\small Chengdu, 610054, P.R. China}\\
{\small E-mail: lvhz08@lzu.edu.cn}\\
{\small $^{2}$School of Mathematics and Information Science, North Minzu University, }\\
{\small Yinchuan, 750027, P.R. China} \\
{\small E-mail: dengkai04@126.com}\\
{\small $^{3}$School of Mathematics, Southwest Jiaotong University, }\\
{\small Chengdu, 610031, P.R. China} \\
{\small E-mail: yangxiaomath@swjtu.edu.cn}\\}

\date{}
\begin{document}

\maketitle
\begin{abstract}

Li et al. in [Inf. Process. Lett. 77 (2001) 35--41] proposed the shuffle cube $SQ_{n}$ as an attractive interconnection network topology for massive parallel and distributed systems. Symmetric property, including vertex-transitivity and edge-transitivity, are the most desirable property for networks since their routing algorithms are universal for each vertex. By far, vertex-transitivity and edge-transitivity of the shuffle cube remain unknown. In this paper, we show that $SQ_{n}$ is not vertex-transitive for all $n>2$, which is not an appealing property in interconnection networks. To overcome this limitation, two novel vertex-transitive variants of the shuffle-cube, namely simplified shuffle-cube $SSQ_{n}$ and balanced shuffle cube $BSQ_{n}$ are introduced. Then, routing algorithms of $SSQ_{n}$ and $BSQ_{n}$ for all $n>2$ are given respectively. Furthermore, we show that both $SSQ_{n}$ and $BSQ_{n}$ possess Hamiltonian cycle embedding for all $n>2$. Finally, as a by-product, we mend a flaw in the Property 3 in [IEEE Trans. Comput. 46 (1997) 484--490].

\vskip 0.1 in

\noindent \textbf{Key words:} Interconnection network; Shuffle-cube; Vertex-transitive; Simplified shuffle-cube; Balanced shuffle-cube

\noindent \textbf{Mathematics Subject Classification:} 05C60, 68R10
\end{abstract}

\section{Introduction}

In massive processing systems (MPS), processors are connected based on a specific interconnection network topology, which is usually represented by an undirected graph: vertices represent processors and edges represent links between processors. It is well-known that the network topology plays a crucial role in its performance. The well-known hypercube (or $n$-cube), which has received much attention over past decades, is one of the most popular topologies for MPS. The hypercube possesses numerous desirable properties for MPS, such as large bandwidth, short (logarithmic) diameter, high symmetry and connectivity, recursive structure, simple routing and broadcasting. With so many advantages, a number of hypercube machines have been implemented \cite{Hayes}, such as Cosmic Cube \cite{Seitz}, the Ametek S/14 \cite{Athas}, the iPSC \cite{Dunican}, the NCUBE \cite{Ncube}, and the CM-200 \cite{Brunet}.

However, a network topology can not be optimum in all aspects. To enhance some properties of the hypercube, a number of variations of the hypercube have been proposed. One method of reducing diameter of the hypercube is to still retain the ``hypercube-like'' structure by twisting some pair of edges, such as twisted cube \cite{Abraham}, M\"{o}bius cube \cite{Cull}, crossed cube \cite{Huang} and Mcube \cite{Singhvi}. Another feasible approach is to reduce vertex degree that leads to technological problems in parallel computing, such as Fibonacci cube \cite{Hsu}, Lucas cube \cite{Munarini} and exchanged hypercube \cite{Loh}. In addition, other pleasing hypercube variants are obtained by adopting different techniques, such as balanced hypercube \cite{Wu}, folded hypercube \cite{El-Amawy}, generalized hypercube \cite{Bhuyan} and shuffle-cube \cite{Li}. The shuffle-cube $SQ_n$ has some good combinatorial properties and fault-tolerant properties. The reliability of the shuffle cube concerning refined connectivity were determined in \cite{Xu,Qin,Ding}. The conditional diagnosability of $SQ_n$ was studied in \cite{Xu2,Lin}. The matching preclusion number was determined by Antantapantula et al. \cite{Antantapantula}. Fault-tolerant Hamiltonian cycle embedding of the shuffle-cube was investigated by Li et al. \cite{Li2}.

To achieve high performance, graphs with high levels of symmetry (e.g., Cayley graphs \cite{Akers,Lakshmivarahan}) are recommended as network topologies since it often simplifies the computation and routing algorithms. Actually, vertex-transitive and/or edge-transitive graphs are widely used to design networks of high levels of symmetry. For example, numerous attractive networks, including hypercube, $k$-ary $n$-cube, balanced hypercube, folded hypercube, star graphs and some of their variants, which have both theoretical and practical importance, are vertex-transitive and/or edge-transitive. On the other hand, to study the symmetry of a graph, the aim is to obtain as much information as possible about its symmetric property.

{\em Routing} is the problem of finding path between each pair of vertices in a graph. It is well-known that a good routing algorithm does not require large amount of memory resources to build the routing table, and its convergence time is usually slow for large networks. Because of the importance of a good routing in interconnection networks, it has received considerable attention in the literature \cite{Libeskind-Hadas,Liu,Khosravi,Yu}. In particular, a path/cycle consisting of all nodes of a network, i.e. Hamiltonian path/cycle, has a wide range of importance in theory and practice since there are some applications of Hamiltonian path in the on-line optimization of a complex Flexible Manufacturing System \cite{Ascheuer}, as well as full utilization of all nodes in a network \cite{Park}.

The rest of this paper is organized as follows. In Section 2, some notations and the definition of the shuffle-cube are presented. In Section 3, it is showed that the shuffle-cube is neither vertex-transitive nor edge-transitive. In Section 4, two vertex-transitive variants of the shuffle cube are proposed. The routing algorithms and Hamiltonian cycle embeddings of the variants are investigated in Section 5 and 6, respectively. Conclusions are given in Section 7.

\section{Preliminaries}

Let $G=(V(G),E(G))$ be a graph, where $V(G)$ is the vertex-set of $G$ and $E(G)$ is the edge-set of $G$. The number of vertices of $G$ is denoted by $|G|$. A path $P=\langle x_0,x_1,\cdots, x_k\rangle$ in $G$ is a sequence of distinct vertices so that there is an edge joining each pair of consecutive vertices. If a path $C=\langle x_0,x_1,\cdots, x_k\rangle$ is such that $k\geq3$, $x_0=x_k$, then $C$ is said to be a {\em cycle}, and the length of $C$ is the number of edges contained in $C$. The length of the shortest cycle of $G$ is called the {\em girth} of $G$, denoted by $g(G)$. In particular, a cycle containing all vertices of $G$ is called a {\em Hamiltonian cycle}. The {\em clique} of $G$ is a set of pairwise adjacent vertices and the {\em clique number} of $G$, is the maximum size of a clique in $G$. A graph $G$ is vertex-transitive if for each pair $u,v\in V(G)$ there exists an automorphism that maps $u$ to $v$. A graph $G$ is edge-transitive if for all $e,f\in E(G)$ there
exists an automorphism of $G$ that maps the endpoints of $e$ to the endpoints of $f$. For other standard graph notations not defined here please refer to \cite{Bondy}.

The vertices are labelled by binary sequences of $n$-bit. For a vertex $u=u_{n-1}u_{n-2}$ $\cdots u_1u_0$, $u_i\in\{0,1\}$ for each $0\leq i\leq n-1$, the $j$-{\em prefix} of $u$ is $u_{n-1}u_{n-2}\cdots$ $u_{n-j}$, written by $p_j(u)$, and the $k$-{\em suffix} of $u$ is $u_{k-1}u_{k-2}\cdots u_1u_0$, written by $s_k(u)$. The {\em Hamming distance} of two vertices $u$ and $v$, denoted by $h(u,v)$, is the number of bits which they differ. The well-known $n$-dimensional hypercube $Q_n$, consists of all of the $n$-bit binary sequences as its vertex set and two distinct vertices $u$ and $v$ are linked by an edge if and only if $h(u,v)=1$.

To recursively build shuffle-cubes, we define four sets containing $4$-tuple of binary sequences as follows:

\begin{itemize}
\item $V_{00}=\{1111,0001,0010,0011\}$,
\item $V_{01}=\{0100,0101,0110,0111\}$,
\item $V_{10}=\{1000,1001,1010,1011\}$,
\item $V_{11}=\{1100,1101,1110,1111\}$.
\end{itemize}

We are ready to give the definition of the shuffle-cube.
\vskip 0.0 in

\begin{definition}{\bf .}\label{def-shuffle}\cite{Li}
The $n$-{\em dimensional shuffle-cube}, $SQ_n$, is recursively defined as follows: $SQ_2\cong Q_2$. For $n\geq3$, $SQ_n$ contains exactly 16 subcubes $SQ_{n-4}^{i_1i_2i_3i_4}$, where $i_1,i_2,i_3,i_4\in \{0,1\}$ and all vertices of $SQ_{n-4}^{i_1i_2i_3i_4}$ share the same $p_4(u)=i_1i_2i_3i_4$. The vertices $u=u_{n-1}u_{n-2}\cdots u_1u_0$ and $v=v_{n-1}v_{n-2}\cdots v_1v_0$ in different subcubes of dimension $n-4$ are adjacent in $SQ_n$ iff

\begin{enumerate}
\item $s_{n-4}(u)=s_{n-4}(v)$, and
\item $p_{4}(u)\oplus p_{4}(v)\in V_{s_2(u)}$,
\end{enumerate}
where the notation ``$\oplus$'' means bitwise addition under modulo 2.
\end{definition}

\begin{figure}
\centering
\includegraphics[height=80mm]{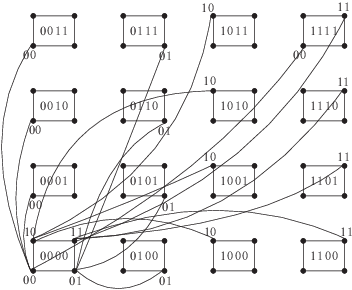}
\caption{$SQ_6$.} \label{SQ6}
\end{figure}

By the definition above, it is clear that $SQ_n$ is $n$-regular and $|V(SQ_n)|=2^n$, which is the same as that of $Q_n$. For clarity, $SQ_{6}$ is illustrated in Fig. \ref{SQ6} with only edges incident to vertices in $SQ_2^{0000}$.

By Definition \ref{def-shuffle}, it implies that two vertices in different subcubes of dimension $n-4$ differ in exactly one 4-bit. For convenience, the $j$-th 4-bit of a vertex $u$, denoted by $u_4^j$, is defined as $u_4^j=u_{4j+1}u_{4j}u_{4j-1}u_{4j-2}$, $1\leq j\leq \frac{n-2}{4}$. For notation consistency, we define $u_4^0=u_{1}u_{0}$. For two distinct vertices $u$ and $v$ in $SQ_n$, $u_4^{j}=v_4^{j}$ if they bitwise equal. The complementary of an arbitrary bit $u_i$ of $u$, i.e. $1-u_i$, is denoted by $\overline{u}_i$. Similar to Hamming distance, we define 4-bit Hamming distance between $u$ and $v$, written by $h_4(u,v)$, as the number of 4-bits $u^j_4$ with $0\leq j\leq \frac{n-2}{4}$ such that $u^j_4\neq v^j_4$. In particular, we use $h^*_4(u,v)$ to denote the number of 4-bits $u^j_4$ with $1\leq j\leq \frac{n-2}{4}$ such that $u^j_4\neq v^j_4$.

\section{Symmetric properties of shuffle-cubes}
In this section, we shall state some symmetric properties of shuffle-cubes. We begin with the following lemma.

\begin{lemma}\label{bipartite}{\bf.} $SQ_{n}$ is non-bipartite for all $n>2$. Moreover, $g(SQ_n)=3$ for all $n>2$.
\end{lemma}
\noindent{\bf Proof.} To show that $SQ_{n}$ is non-bipartite for all $n>2$, it suffices to present an odd cycle of $SQ_n$. Let $u=u_{n-1}u_{n-2}\cdots u_1u_0$ be a vertex in $SQ_{n}$ with $u_4^0=u_1u_0=00$ and let $v,w,x$ and $y$ be four neighbors of $u$ with $u_4^j\neq v_4^j, w_4^j, x_4^j, y_4^j$ for some $1\leq j\leq\frac{n-2}{4}$. It follows that $v_4^0=w_4^0=x_4^0=y_4^0=00$.

Additionally, we may assume that $u_{4j+1}u_{4j}=v_{4j+1}v_{4j}=w_{4j+1}w_{4j}=x_{4j+1}x_{4j}$ and $y_{4j+1}y_{4j}=\bar{u}_{4j+1}\bar{u}_{4j}$. It follows from Definition \ref{def-shuffle} that $u,v,w$ and $x$ are pairwise adjacent, forming a clique of size four. Moreover, $y$ is not adjacent to none of $v,w$ and $x$.
Suppose that there exists another vertex $z$ ($\neq y$) such that $z$ is adjacent to all of $u,v,w$ and $x$. We claim that $z_4^j\neq u_4^j,v_4^j,w_4^j,x_4^j$. Otherwise, suppose without loss of generality that $z_4^j=u_4^j$. Therefore, $z_4^i\neq u_4^i$ for some $i\neq j$, $0\leq i\leq \frac{n-2}{4}$. So $z$ and $v$ (resp. $w$, $x$) differ in two 4-bits, a contradiction. Thus, the claim holds. Note that $u_4^j\oplus z_4^j\in V_{00}$, then $z_4^j=\overline{u}_{4j+1}\overline{u}_{4j}\overline{u}_{4j-1}\overline{u}_{4j-2}$. By Definition \ref{def-shuffle}, $z$ is not adjacent to $v,w$ and $x$, a contradiction again.

This completes the proof. \qed

\begin{lemma}\label{K-4}{\bf.} A vertex of $u=u_{n-1}u_{n-2}\cdots u_1u_0$ of $SQ_{n}$ ($n>2$) is contained in a clique of size four if and only if $u_4^0=00$. Moreover, the clique number of $SQ_{n}$ is four.
\end{lemma}
\noindent{\bf Proof.} It has been shown in the proof of Lemma \ref{bipartite} that if $u_4^0=00$, there are four vertices $u,v,w$ and $x$ forming a clique of size four. To prove the clique number of is four, we show that there exists no vertex adjacent to all of $u,v,w$ and $x$.

In what follows, we prove the necessity. Suppose on the contrary that $u_4^0\neq00$. We shall prove that any pair of the neighbors of $u$ are nonadjacent. Let $v$ and $w$ be two arbitrary neighbors of $u$ and let $u_4^i\neq v_4^i$ and $u_4^j\neq w_4^j$, where $0\leq i,j\leq\frac{n-2}{4}$. The following two cases arise.

\noindent{\bf Case 1:} $i\neq j$. This implies that $v$ and $w$ differ in two $4$-bits. By Definition \ref{def-shuffle}, $v$ and $w$ are nonadjacent.

\noindent{\bf Case 2:} $i=j$. If $i=j=0$, then $v$ and $w$ are in the same $SQ_2$ of $SQ_n$. By Definition \ref{def-shuffle}, $v$ and $w$ are nonadjacent. So we assume that $i=j>0$. Clearly, $u_4^0=v_4^0=w_4^0$. It follows that $u_{4i+1}u_{4i}\oplus v_{4i+1}v_{4i}=u_4^0$ and $u_{4i+1}u_{4i}\oplus w_{4i+1}w_{4i}=u_4^0$, implying that $v_{4i+1}v_{4i}=w_{4i+1}w_{4i}\neq u_4^0$. On the contrary, suppose that $v$ and $w$ are adjacent. Since $v_4^i\neq w_4^i$, we have $v_4^i\oplus w_4^i=00xy$, where $x,y\in\{0,1\}$ are two bits we do not care. This contradicts the fact that the first two bits of $v_4^i$, $w_4^i$ are the same as $v_4^0$. This completes the proof. \qed

\vskip 0.1 in

\noindent{\bf Remark 1.}\label{4-clique} Any vertex $u=u_{n-1}u_{n-2}\cdots u_1u_0$ of $SQ_n$ ($n>2$) with $u_1u_0=00$ and its $j$-bit neighbors $v,w$ and $x$ with $u_{4j+1}u_{4j}=v_{4j+1}v_{4j}=w_{4j+1}w_{4j}=x_{4j+1}x_{4j}$, form a clique of size four. So $u$ is contained in exactly $\frac{n-2}{4}$ cliques of size four. Moreover, $SQ_n$ consists of $4^{\frac{n-2}{4}}$ vertex-disjoint clique of size four.

\begin{theorem}{\bf.}\label{not-v-t} $SQ_{n}$ is vertex-transitive for $n=2$, but not vertex-transitive for all $n>2$.
\end{theorem}
\noindent{\bf Proof.} Clearly, $SQ_{2}$ consists of 16 vertex-disjoint 4-cycles, which is obviously vertex-transitive. By Lemma \ref{K-4}, we know that each vertex $u$ with $u_4^0=00$ is contained in a clique of size four in $SQ_n$, $n>2$. Any vertex $v$ with $u_4^0\neq00$ is contained in no clique of size four, showing obviously that $SQ_{n}$ is not vertex-transitive for all $n>2$. \qed

It is well-known that an edge-transitive but not vertex-transitive graph must be bipartite. Combing Lemma \ref{bipartite} with Theorem \ref{not-v-t}, we have the following corollary.

\begin{corollary}{\bf.}\label{not-e-t} $SQ_{n}$ is not edge-transitive for all $n>2$.
\end{corollary}

\section{Two vertex-transitive variations of the shuffle-cubes}

As we have shown in the previous section, the shuffle-cube $SQ_{n}$ ($n>2$) is not vertex-transitive. In this section, we shall present two variants of the shuffle-cubes, which are both vertex-transitive. We shall adopt the notations which we have defined in the $SQ_{n}$.

\begin{definition}{\bf.}\label{variation-1} The $n$-{\em dimensional simplified shuffle cube}, $SSQ_n$, is recursively defined as follows: $SSQ_2\cong Q_2$. For $n>2$, $SSQ_n$ contains exactly 8 subcubes $SSQ_{n-4}^{i_1i_2i_3i_4}$, where $i_1i_2\in\{00,11\}$, $i_3,i_4\in \{0,1\}$ and all vertices of $SSQ_{n-4}^{i_1i_2i_3i_4}$ share the same $p_4(u)=i_1i_2i_3i_4$. The vertices $u=u_{n-1}u_{n-2}\cdots u_1u_0$ and $v=v_{n-1}v_{n-2}\cdots v_1v_0$ in different subcubes of dimension $n-4$ are adjacent in $SSQ_n$ iff

\begin{enumerate}
\item $s_{n-4}(u)=s_{n-4}(v)$, and
\item $p_{4}(u)\oplus p_{4}(v)\in V_{00}$,
\end{enumerate}
where the notation ``$\oplus$'' means bitwise addition under modulo 2.
\end{definition}

\begin{figure}
\centering
\includegraphics[height=80mm]{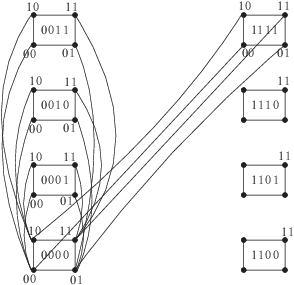}
\caption{$SSQ_6$.} \label{SSQ6}
\end{figure}

The reason why we define $i_1i_2\in\{00,11\}$ (which is slightly different from the definition of $SQ_n$) is that if we define $i_1,i_2 \in\{0,1\}$ as in Definition \ref{def-shuffle}, then the resulting graph is always disconnected, which does not meet the basic requirement of interconnection networks. As a result, $SSQ_n$ has $2^{\frac{3n+2}{4}}$ vertices, $n\equiv 2(\text{mod}\ 4)$. Clearly, $SSQ_n$ is also $n$-regular and non-bipartite. $SSQ_{6}$ is illustrated in Fig. \ref{SSQ6} with only edges incident to vertices in $SSQ_2^{0000}$.

%

In what follows, we shall prove the vertex-transitivity of $SSQ_n$, $n\geq2$.

\begin{theorem}{\bf.}
$SSQ_n$ is vertex-transitive whenever $n\geq2$.
\end{theorem}
\noindent{\bf Proof.} Obviously, $SSQ_n$ is vertex-transitive if $n=2$. It suffices to show that $SSQ_n$ is vertex-transitive for $n>2$.

Let $u=u_{n-1}u_{n-2}\cdots u_1u_0$ and $v=v_{n-1}v_{n-2}\cdots v_1v_0$ be any two distinct vertices of $SSQ_n$. Our objective is to show that there is an automorphism $\phi$ of $SSQ_n$ such that $\phi(v)=u$. Thus, for any vertex $w=w_{n-1}w_{n-2}\cdots w_1w_0$ of $SSQ_n$, we define $\phi$ 4-bit by 4-bit (in fact, $w_4^0$ contains exactly two bits) below:

$$\phi(w_4^j)=w_4^j\oplus(u_4^j\oplus v_4^j),\ 0\leq j\leq\frac{n-2}{4}.$$

Now we are ready to show that $\phi$ is an isomorphism of $SSQ_n$. Let $xy$ be an arbitrary edge of $SSQ_n$. We shall show that $\phi(x)\phi(y)$ is also an edge of $SSQ_n$. By Definition \ref{variation-1}, we known that $x=x_{n-1}x_{n-2}\cdots x_1x_0$ and $y=y_{n-1}y_{n-2}\cdots y_1y_0$ differ in exactly one 4-bit (including the last two bits), say $j$, $0\leq j\leq\frac{n-2}{4}$.

If $x_4^0\neq y_4^0$, then $x_1x_0\neq y_1y_0$ and $x_{n-1}x_{n-2}\cdots x_{3}=y_{n-1}y_{n-2}\cdots y_{3}$. By the definition of $\phi$, we have $\phi(x_4^0)=x_4^0\oplus(u_4^j\oplus v_4^j)$ and $\phi(y_4^0)=y_4^0\oplus(u_4^j\oplus v_4^j)$. Additionally, $\phi(x_4^j)=\phi(y_4^j)$ for $1\leq j\leq\frac{n-2}{4}$. Thus, $\phi(x)\phi(y)\in E(SSQ_n)$.

If $x_4^j\neq y_4^j$ for some $j\in\{1,\cdots,\frac{n-2}{4}\}$, then $x_1x_0=y_1y_0$ and $x_4^i=y_4^i$ for all $i\in\{0,\cdots,\frac{n-2}{4}\}\setminus\{j\}$.

By the definition of $\phi$, we have $\phi(x_4^j)=x_4^j\oplus(u_4^j\oplus v_4^j)$ and $\phi(y_4^j)=y_4^j\oplus(u_4^j\oplus v_4^j)$. Additionally, $\phi(x_4^i)=\phi(y_4^i)$ for all $i\in\{0,\cdots,\frac{n-2}{4}\}\setminus\{j\}$. Thus, $\phi(x)\phi(y)\in E(SSQ_n)$. This completes the proof.\qed

In the following, we shall give another variation of the shuffle cubes, which is motivated by the idea of construction of the balanced hypercubes \cite{Wu}. For the readability of this paper, we first present the definition of the balanced hypercube, and then give the definition of the balanced shuffle cube.

\begin{definition}\cite{Wu}\label{def-balanced-cube}{\bf.}
An $n$-dimensional balanced hypercube $BH_{n}$ contains
$4^{n}$ vertices $(a_{0},$ $\cdots,a_{i-1},$
$a_{i},a_{i+1},\cdots,a_{n-1})$, where $a_{i}\in\{0,1,2,3\}$ $(0\leq
i\leq n-1)$. Any vertex $v=(a_{0},\cdots,a_{i-1},$
$a_{i},a_{i+1},\cdots,a_{n-1})$ in $BH_{n}$ has the following $2n$ neighbors:

\begin{enumerate}
\item $((a_{0}+1)$ mod $
4,a_{1},\cdots,a_{i-1},a_{i},a_{i+1},\cdots,a_{n-1})$,\\
      $((a_{0}-1)$ mod $ 4,a_{1},\cdots,a_{i-1},a_{i},a_{i+1},\cdots,a_{n-1})$, and
\item $((a_{0}+1)$ mod $ 4,a_{1},\cdots,a_{i-1},(a_{i}+(-1)^{a_{0}})$ mod $
4,a_{i+1},\cdots,a_{n-1})$,\\
      $((a_{0}-1)$ mod $ 4,a_{1},\cdots,a_{i-1},(a_{i}+(-1)^{a_{0}})$ mod $
      4,a_{i+1},\cdots,a_{n-1})$.
\end{enumerate}
\end{definition}

\begin{definition}{\bf.}\label{variation-2} The $n$-{\em dimensional balanced shuffle cube}, $BSQ_n$, is recursively defined as follows: $BSQ_2\cong Q_2$. For $n\geq3$, $BSQ_n$ contains exactly 16 subcubes $BSQ_{n-4}^{i_1i_2i_3i_4}$, where $i_1,i_2,i_3,i_4\in \{0,1\}$ and all vertices of $BSQ_{n-4}^{i_1i_2i_3i_4}$ share the same $p_4(u)=i_1i_2i_3i_4$. The vertices $u=u_{n-1}u_{n-2}\cdots u_1u_0$ and $v=v_{n-1}v_{n-2}\cdots v_1v_0$ in different subcubes of dimension $n-4$ are adjacent in $BSQ_n$ iff

\begin{enumerate}
\item $s_{n-4}(u)=s_{n-4}(v)$, and
\item $u_{n-2}$ and $v_{n-2}$ have different parities, and $v_{n-1}v_{n-2}=u_{n-1}u_{n-2}\pm1$, and
\item $v_{n-3}v_{n-4}=u_{n-3}u_{n-4}+(-1)^{u_{n-2}}$, or $v_{n-3}v_{n-4}=u_{n-3}u_{n-4}$.
\end{enumerate}
where addition and substraction are under modulo 4 by regarding the two bits as an integer.
\end{definition}

\begin{figure}
\centering
\includegraphics[height=80mm]{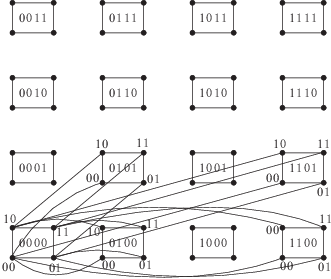}
\caption{$BSQ_6$.} \label{BSQ6}
\end{figure}

Clearly, $BSQ_n$ has $2^{n}$ vertices, $n\equiv 2(\text{mod}\ 4)$. Additionally, $BSQ_n$ is $n$-regular and bipartite. $BSQ_{6}$ is illustrated in Fig. \ref{BSQ6} with only edges incident to vertices in $BSQ_2^{0000}$.

We shall prove the following statement, which directly leads to the regularity of $BSQ_n$, $n\geq2$.

\begin{theorem}{\bf.}\label{Cayley-BSQn}
$BSQ_n$ is vertex-transitive whenever $n>2$.
\end{theorem}
\noindent{\bf Proof.} Clearly, $BSQ_n$ is vertex-transitive when $n=2$. It remains to show that $BSQ_n$ is vertex-transitive for $n>2$.

Let $v=v_{n-1}v_{n-2}\cdots v_1v_0$ and $u=u_{n-1}u_{n-2}\cdots u_1u_0$ be any two distinct vertices of $BSQ_n$. Our aim is to show that there is an automorphism $\psi$ of $BSQ_n$ such that $\psi(v)=u$. Thus, for any vertex $w=w_{n-1}w_{n-2}\cdots w_1w_0$ of $BSQ_n$, we define $\psi$ 4-bit by 4-bit (in fact, $w_4^0$ contains exactly two bits) below:\\
If both of $u_{4j}$ and $v_{4j}$ are even, then
$$\psi(w_4^j)=(v_{4j+1}v_{4j}-u_{4j+1}u_{4j}+w_{4j+1}w_{4j},v_{4j-1}v_{4j-2}-u_{4j-1}u_{4j-2}+w_{4j-1}w_{4j-2});$$
If $u_{4j}$ is even and $v_{4j}$ is odd, then
$$\psi(w_4^j)=(v_{4j+1}v_{4j}-u_{4j+1}u_{4j}-w_{4j+1}w_{4j},v_{4j-1}v_{4j-2}-u_{4j-1}u_{4j-2}-w_{4j-1}w_{4j-2});$$
If $u_{4j}$ is odd and $v_{4j}$ is even, then
$$\psi(w_4^j)=(-v_{4j+1}v_{4j}+u_{4j+1}u_{4j}-w_{4j+1}w_{4j},-v_{4j-1}v_{4j-2}+u_{4j-1}u_{4j-2}-w_{4j-1}w_{4j-2});$$
If both of $u_{4j}$ and $v_{4j}$ are odd, then
$$\psi(w_4^j)=(v_{4j+1}v_{4j}-u_{4j+1}u_{4j}+w_{4j+1}w_{4j},v_{4j-1}v_{4j-2}-u_{4j-1}u_{4j-2}+w_{4j-1}w_{4j-2});$$
$1\leq j\leq\frac{n-2}{4}$, and addition and substraction are under modulo 4 by regarding the two bits as an integer. In addition,
$$\psi(w_4^0)=v_{1}v_{0}+(w_{1}w_{0}-u_{1}u_{0})\ \text{when}\ j=0.$$

Let $xy$ be an arbitrary edge of $BSQ_n$ with $x=x_{n-1}x_{n-2}\cdots x_1x_0$ and $y=y_{n-1}y_{n-2}\cdots y_1y_0$. We shall show that $\psi(x)\psi(y)$ is also an edge of $BSQ_n$. By Definition \ref{variation-2}, it is obvious that $x$ and $y$ differ in exactly one 4-bit (including the last two bits), say $j$, $0\leq j\leq\frac{n-2}{4}$.

If $x_4^0\neq y_4^0$, then $x_1x_0\neq y_1y_0$ and $x_{n-1}x_{n-2}\cdots x_{3}=y_{n-1}y_{n-2}\cdots y_{3}$. By the definition of $\psi$, we have $\psi(x_4^0)=v_{1}v_{0}+(x_{1}x_{0}-u_{1}u_{0})$ and $\phi(y_4^0)=v_{1}v_{0}+(y_{1}y_{0}-u_{1}u_{0})$. Additionally, $\psi(x_4^j)=\psi(y_4^j)$ for $1\leq j\leq\frac{n-2}{4}$. Thus, $\psi(x)\psi(y)\in E(BSQ_n)$.

If $x_4^j\neq y_4^j$ for some $j\in\{1,\cdots,\frac{n-2}{4}\}$, then $x_1x_0=y_1y_0$ and $x_4^i=y_4^i$ for all $i\in\{0,\cdots,\frac{n-2}{4}\}\setminus\{j\}$. By the definition of $\psi$, the following cases arise.

\noindent{\bf Case 1.} $u_{4j}$ is even.

\noindent{\bf Case 1.1.} $v_{4j}$ is even. Clearly, we have\\
$\psi(x_4^j)=(v_{4j+1}v_{4j}-u_{4j+1}u_{4j}+x_{4j+1}x_{4j},v_{4j-1}v_{4j-2}-u_{4j-1}u_{4j-2}+x_{4j-1}x_{4j-2})$
and\\
$\psi(y_4^j)=(v_{4j+1}v_{4j}-u_{4j+1}u_{4j}+y_{4j+1}y_{4j},v_{4j-1}v_{4j-2}-u_{4j-1}u_{4j-2}+y_{4j-1}y_{4j-2})$. Thus, the parity of $v_{4j+1}v_{4j}-u_{4j+1}u_{4j}+x_{4j+1}x_{4j}$ (resp. $v_{4j+1}v_{4j}-u_{4j+1}u_{4j}+y_{4j+1}y_{4j}$) is the same as that of $x_{4j+1}x_{4j}$ (resp. $y_{4j+1}y_{4j}$). It follows from Definition \ref{variation-2} that $\psi(x)\psi(y)\in E(BSQ_n)$.

\noindent{\bf Case 1.2.} $v_{4j}$ is odd. Clearly, we have\\ $\psi(x_4^j)=(v_{4j+1}v_{4j}-u_{4j+1}u_{4j}-x_{4j+1}x_{4j},v_{4j-1}v_{4j-2}-u_{4j-1}u_{4j-2}-x_{4j-1}x_{4j-2})$ and
$\psi(y_4^j)=(v_{4j+1}v_{4j}-u_{4j+1}u_{4j}-y_{4j+1}y_{4j},v_{4j-1}v_{4j-2}-u_{4j-1}u_{4j-2}-y_{4j-1}y_{4j-2})$. Observe that $xy\in E(BSQ_n)$.

If $x_{4j}$ is even, then $y_{4j-1}y_{4j-2}\equiv x_{4j-1}x_{4j-2}+1(\text{mod}\ 4)$. Thus, $-y_{4j-1}y_{4j-2}\equiv -x_{4j-1}x_{4j-2}-1(\text{mod}\ 4)$. So $v_{4j-1}v_{4j-2}-u_{4j-1}u_{4j-2}-y_{4j-1}y_{4j-2}\equiv v_{4j-1}v_{4j-2}-u_{4j-1}u_{4j-2}-x_{4j-1}x_{4j-2}-1(\text{mod}\ 4)$. Now we see that $v_{4j+1}v_{4j}-u_{4j+1}u_{4j}-x_{4j+1}x_{4j}$ is odd and $v_{4j+1}v_{4j}-u_{4j+1}u_{4j}-y_{4j+1}y_{4j}$ is even. By Definition \ref{variation-2}, we have $\psi(x)\psi(y)\in E(BSQ_n)$.

If $x_{4j}$ is odd, then $y_{4j-1}y_{4j-2}\equiv x_{4j-1}x_{4j-2}-1(\text{mod}\ 4)$. Thus, $-y_{4j-1}y_{4j-2}\equiv -x_{4j-1}x_{4j-2}+1(\text{mod}\ 4)$. So $v_{4j-1}v_{4j-2}-u_{4j-1}u_{4j-2}-y_{4j-1}y_{4j-2}\equiv v_{4j-1}v_{4j-2}-u_{4j-1}u_{4j-2}-x_{4j-1}x_{4j-2}+1(\text{mod}\ 4)$. Now we see that $-v_{4j+1}v_{4j}+u_{4j+1}u_{4j}-x_{4j+1}x_{4j}$ is even and $-v_{4j+1}v_{4j}+u_{4j+1}u_{4j}-y_{4j+1}y_{4j}$ is odd. By Definition \ref{variation-2}, we have $\psi(x)\psi(y)\in E(BSQ_n)$.

\noindent{\bf Case 2.} $u_{4j}$ is odd.

\noindent{\bf Case 2.1.} $v_{4j}$ is even. In this case, we have\\
$\psi(x_4^j)\!=\!(\!-v_{4j+1}\!v_{4j}\!+u_{4j+1}\!u_{4j}-x_{4j+1}\!x_{4j},-v_{4j-1}\!v_{4j-2}\!+u_{4j-1}\!u_{4j-2}-x_{4j-1}\!x_{4j-2}),$
and\\
$\psi(y_4^j)\!=\!(\!-v_{4j+1}\!v_{4j}\!+u_{4j+1}\!u_{4j}-y_{4j+1}\!y_{4j},-v_{4j-1}\!v_{4j-2}\!+u_{4j-1}\!u_{4j-2}-y_{4j-1}\!y_{4j-2})$. Obviously, $-v_{4j+1}v_{4j}+u_{4j+1}u_{4j}-x_{4j+1}x_{4j}$ and $-v_{4j+1}v_{4j}+u_{4j+1}u_{4j}-y_{4j+1}y_{4j}$ have different parities. Observe that $xy\in E(BSQ_n)$.

If $x_{4j}$ is even, then $y_{4j-1}y_{4j-2}\equiv x_{4j-1}x_{4j-2}+1(\text{mod}\ 4)$. Thus, $-y_{4j-1}y_{4j-2}\equiv -x_{4j-1}x_{4j-2}-1(\text{mod}\ 4)$. Now we see that $-v_{4j+1}v_{4j}+u_{4j+1}u_{4j}-x_{4j+1}x_{4j}$ is odd and $-v_{4j+1}v_{4j}+u_{4j+1}u_{4j}-y_{4j+1}y_{4j}$ is even. By Definition \ref{variation-2}, we have $\psi(x)\psi(y)\in E(BSQ_n)$.

If $x_{4j}$ is odd, then $y_{4j-1}y_{4j-2}\equiv x_{4j-1}x_{4j-2}-1(\text{mod}\ 4)$. Thus, $-y_{4j-1}y_{4j-2}\equiv -x_{4j-1}x_{4j-2}+1(\text{mod}\ 4)$. Now we see that $-v_{4j+1}v_{4j}+u_{4j+1}u_{4j}-x_{4j+1}x_{4j}$ is even and $-v_{4j+1}v_{4j}+u_{4j+1}u_{4j}-y_{4j+1}y_{4j}$ is odd. By Definition \ref{variation-2}, we have $\psi(x)\psi(y)\in E(BSQ_n)$.

\noindent{\bf Case 2.2.} $v_{4j}$ is odd. In this case, we have\\
$\psi(x_4^j)=(v_{4j+1}v_{4j}-u_{4j+1}u_{4j}+x_{4j+1}x_{4j},v_{4j-1}v_{4j-2}-u_{4j-1}u_{4j-2}+x_{4j-1}x_{4j-2})$, and\\
$\psi(y_4^j)=(v_{4j+1}v_{4j}-u_{4j+1}u_{4j}+y_{4j+1}y_{4j},v_{4j-1}v_{4j-2}-u_{4j-1}u_{4j-2}+y_{4j-1}y_{4j-2})$.

If $x_{4j}$ is even, then $y_{4j-1}y_{4j-2}\equiv x_{4j-1}x_{4j-2}+1(\text{mod}\ 4)$. We see that $-v_{4j+1}v_{4j}+u_{4j+1}u_{4j}-x_{4j+1}x_{4j}$ is even and $-v_{4j+1}v_{4j}+u_{4j+1}u_{4j}-y_{4j+1}y_{4j}$ is odd. By Definition \ref{variation-2}, we have $\psi(x)\psi(y)\in E(BSQ_n)$.

If $x_{4j}$ is odd, then $y_{4j-1}y_{4j-2}\equiv x_{4j-1}x_{4j-2}-1(\text{mod}\ 4)$. Now we see that $-v_{4j+1}v_{4j}+u_{4j+1}u_{4j}-x_{4j+1}x_{4j}$ is odd and $-v_{4j+1}v_{4j}+u_{4j+1}u_{4j}-y_{4j+1}y_{4j}$ is even. By Definition \ref{variation-2}, we have $\psi(x)\psi(y)\in E(BSQ_n)$.

This completes the proof.\qed

\section{Routing algorithms of SSQ$_n$ and BSQ$_n$}

In this section, we shall present the following shortest path algorithms of SSQ$_n$ and BSQ$_n$, respectively. We begin with the following assertion.

\begin{theorem}{\bf.}\label{diamter-SSQn} $d(SSQ_n)=2$ for $n=2$ and $d(SSQ_n)=\frac{n-2}{2}+2$ for $n>2$.
\end{theorem}
\noindent{\bf Proof.} It is clear that the statement is true for $n=2$. In what follows, we assume that $n>2$.

Let $u$ be a vertex of $SSQ_n$ with all bits zero. Since $SSQ_n$ is vertex transitive, it suffices to determine $\max d(u,v)$ for any vertex $v\in V(SSQ_n)$. Thus, by Definition \ref{variation-1}, for each $0\leq j\leq\frac{n-2}{4}$, there are at most two steps from $u_4^j$ to $v_4^j$. Therefore, $\max d(u,v)\leq \frac{n-2}{2}+2$.

On the other hand, let $v=\underbrace{1101}\underbrace{1101}\cdots\underbrace{1101}11$ be a vertex of $SSQ_n$. Thus, by Definition \ref{variation-1}, for each $0\leq j\leq\frac{n-2}{4}$, there are exactly two steps from $u_4^j$ to $v_4^j$, and each $u_4^j$ (resp. $v_4^j$) is independent to others. Thus, $d(u,v)=\frac{n-2}{2}+2$.

This completes the proof.\qed

By the above statement, it is straightforward to obtain the following algorithm.

\noindent\rule{\textwidth}{0.5pt}
\noindent {\bf Algorithm A} ({\em Routing algorithm of $SSQ_n$})\\
\noindent {\em Input:} $SSQ_n$ with $n>2$.\\
\noindent {\em Output:} a shortest path from an arbitrary vertex $v$ to $u$ in $SSQ_n$.
\noindent {\em Initialization:} $s\leftarrow h_4^*(v,u)$, $S=\{j|v_4^j\neq u_4^j\ \text{and}\ j>0 \}$
\vskip 0.2 cm
\indent{\bf if} $v=u$ \\
\indent\indent route into $u$    \\
\indent{\bf else } \\
\indent\indent {\bf while} $s>0$ {\bf do}\\
\indent\indent arbitrarily choose $j\in S$\\
\indent\indent\indent{\bf if} $v_{4j+1}v_{4j}=u_{4j+1}u_{4j}$ \\
\indent\indent\indent\indent{\bf then} $v_4^j\leftarrow v_{4j+1}v_{4j}u_{4j-1}u_{4j-2}$ and route into $v$\\
\indent\indent\indent{\bf else}  \\
\indent\indent\indent\indent{\bf if} $v_4^j=\overline{u}_{4j+1}\overline{u}_{4j}\overline{u}_{4j-1}\overline{u}_{4j-2}$\\
\indent\indent\indent\indent\indent{\bf then} $v_4^j\leftarrow \overline{u}_{4j+1}\overline{u}_{4j}\overline{u}_{4j-1}\overline{u}_{4j-2}$ and route into $v$\\
\indent\indent\indent\indent{\bf else} $w\leftarrow v$ \\
\indent\indent\indent\indent $w_4^j\leftarrow \overline{u}_{4j+1}\overline{u}_{4j}\overline{u}_{4j-1}\overline{u}_{4j-2}$ and route into $w$\\
\indent\indent\indent\indent $x\leftarrow w$\\
\indent\indent\indent\indent $x_{4j-1}x_{4j-2}\leftarrow v_{4j-1}v_{4j-2}$ and route into $x$\\
\indent\indent\indent\indent $v\leftarrow x$ \\
\indent\indent $s\leftarrow s-1$\\
\indent\indent $S\leftarrow S\setminus\{j\}$ \\
\indent\indent{\bf end while}\\
\indent\indent{\bf if} $v_1v_0=u_1u_0$\\
\indent\indent\indent route into $u$    \\
\indent\indent{\bf else} route into the neighbor $w$ of $v$ that changes $v_1v_0$ in a cyclic manner with respect to 00, 01, 10, 11 \\
\noindent\rule{\textwidth}{0.5pt}

The correctness of Algorithm A is shown as follows.

\begin{theorem}{\bf.} Algorithm A generates a shortest path from a vertex $u$ to $v$ in $SSQ_n$.
\end{theorem}

\noindent{\bf Proof.} Based on Theorem \ref{diamter-SSQn}, by recursively changing each $u_4^j$ to $v_4^j$ for each $0\leq j\leq\frac{n-2}{4}$, we can obtain a shortest path from $u$ to $v$.\qed

\vskip 0.2in

To determine the diameter of the $BSQ_n$, we have the following result.

\begin{theorem}{\bf.}\label{diamter-BSQn} $d(BSQ_n)=n$ whenever $n\geq2$.
\end{theorem}
\noindent{\bf Proof.} It is clear that the statement is true for $n=2$. In what follows, we assume that $n>2$.

Let $u$ be a vertex of $BSQ_n$ with all bits zero. Since $BSQ_n$ is vertex transitive, it suffices to determine $\max d(u,v)$ for any vertex $v\in V(BSQ_n)$. Thus, by Definition \ref{variation-2}, for each $1\leq j\leq\frac{n-2}{4}$, there are at most four steps from $u_4^j$ to $v_4^j$, and there are exactly two steps from $u_4^0$ to $v_4^0$. Therefore, $\max d(u,v)\leq n$.

On the other hand, let $v=\underbrace{1111}\underbrace{1111}\cdots\underbrace{1111}11$ be a vertex of $BSQ_n$. Thus, by Definition \ref{variation-2}, for each $1\leq j\leq\frac{n-2}{4}$, there are exactly four steps from $u_4^j$ to $v_4^j$, and each $u_4^j$ (resp. $v_4^j$) is independent to others. In addition, there are exactly two steps from $u_4^0$ to $v_4^0$. Thus, $d(u,v)=n$.

This completes the proof.\qed

\begin{theorem}{\bf.}\label{neiborhood-BSQn} Let $u$ and $v$ be two distinct vertices of $BSQ_n$ ($n>2$) such that $u_{4j+1}\neq v_{4j+1}$, $u_{4j}u_{4j-1}u_{4j-2}=v_{4j}v_{4j-1}v_{4j-2}$ for some $j\in\{1,\cdots,\frac{n-2}{4}\}$ and $u_4^i=v_4^i$ for each $i\in\{1,\cdots,\frac{n-2}{4}\}\setminus\{j\}$. Then $u$ and $v$ have the same neighborhood.
\end{theorem}
\noindent{\bf Proof.} It can be easily derived from Definition \ref{variation-2}. \qed

\vskip 0.1in

\noindent{\bf Remark 2.} For convenience, if $u$ and $v$ are two distinct vertices satisfy the condition (thus have the same neighborhood) in the above theorem, then we say $u$ is {\em equivalent} to $v$, denote by $u\equiv v$. Without cause of ambiguity, we also denote $u_{4}^j\equiv v_{4}^j$. Moreover, for an arbitrary vertex $u$, there are exactly $\frac{n-2}{4}$ distinct vertices of $BSQ_n$ ($n>2$) that are equivalent to $u$ respectively.

We present the routing algorithm of $BSQ_n$ as follows.

\noindent\rule{\textwidth}{0.5pt}
\noindent {\bf Algorithm B} ({\em Routing algorithm of $BSQ_n$})\\
\noindent {\em Input:} $BSQ_n$ with $n>2$.\\
\noindent {\em Output:} a shortest path from an arbitrary vertex $v$ to $u$ in $BSQ_n$.\\
\noindent {\em Initialization:} $s\leftarrow h_4^*(v,u)$, $S=\{j|v_4^j\neq u_4^j\ \text{and}\ j>0 \}$
\vskip 0.2 cm
\indent{\bf if} $v=u$ \\
\indent\indent route into $u$    \\
\indent{\bf else } \\
\indent\indent {\bf while} $s>0$ {\bf do}\\
\indent\indent arbitrarily choose $j\in S$\\
\indent\indent\indent{\bf if} $v_{4j-1}v_{4j-2}=u_{4j-1}u_{4j-2}$ \\
\indent\indent\indent\indent{\bf if} $v_{4}^j\equiv u_{4}^j$\\
\indent\indent\indent\indent\indent{\bf then} $s\leftarrow s-1$\\
\indent\indent\indent\indent\indent\indent\ \ \  $S\leftarrow S\setminus\{j\}$ \\
\indent\indent\indent\indent{\bf else} \\
\indent\indent\indent\indent\indent{\bf then} route into the neighbor $w$ of $v$ that changes $v_{4j+1}v_{4j}$ in a cyclic manner with respect to 00, 01, 10, 11\\
\indent\indent\indent\indent\indent\indent\ \ \  $s\leftarrow s-1$\\
\indent\indent\indent\indent\indent\indent\ \ \  $S\leftarrow S\setminus\{j\}$ \\
\indent\indent\indent{\bf else}  \\
\indent\indent\indent\indent{\bf if} ($v_{4j-1}v_{4j-2}-u_{4j-1}u_{4j-2}=1\wedge v_{4j}$ is odd)$\vee (v_{4j-1}v_{4j-2}-u_{4j-1}u_{4j-2}=3\wedge v_{4j})$  is even\\
\indent\indent\indent\indent\indent{\bf then} $v_4^j\leftarrow u_4^j$ and route into $v$\\
\indent\indent\indent\indent\indent\indent\ \ \ $s\leftarrow s-1$\\
\indent\indent\indent\indent\indent\indent\ \ \ $S\leftarrow S\setminus\{j\}$ \\
\indent\indent\indent\indent{\bf else if} ($v_{4j-1}v_{4j-2}-u_{4j-1}u_{4j-2}=1\wedge v_{4j}$ is even)$\vee (v_{4j-1}v_{4j-2}-u_{4j-1}u_{4j-2}=3\wedge v_{4j})$ is odd\\
\indent\indent\indent\indent\indent{\bf then} $v_{4j+1}v_{4j}\leftarrow v_{4j+1}v_{4j}+1$ and route into $v$\\
\indent\indent\indent\indent{\bf else}\\
\indent\indent\indent\indent\indent\indent\ \ \ $v_{4j+1}v_{4j}\leftarrow v_{4j+1}v_{4j}+1$ and $v_{4j-1}v_{4j-2}\leftarrow v_{4j-1}v_{4j-2}+(-1)^{v_{4j}}$\\
\indent\indent\indent\indent\indent\indent\ \ \ route into $v$\\
\indent\indent{\bf end while}\\
\indent\indent{\bf if} $v_1v_0=u_1u_0$\\
\indent\indent\indent route into $u$    \\
\indent\indent{\bf else} route into the neighbor $w$ of $v$ that changes $v_1v_0$ in a cyclic manner with respect to 00, 01, 10, 11\\
\noindent\rule{\textwidth}{0.5pt}

\vskip 0.1in

\begin{theorem}{\bf.} Algorithm B generates a shortest path from a vertex $u$ to $v$ in $BSQ_n$.
\end{theorem}

\noindent{\bf Proof.} Based on Theorem \ref{diamter-BSQn}, by recursively changing each $u_4^j$ to $v_4^j$ for each $0\leq j\leq\frac{n-2}{4}$, we can obtain a shortest path from $u$ to $v$.\qed

\section{Hamiltonian cycle embedding of $SSQ_n$ and $BSQ_n$}

The Hamiltonian ring embedding of interconnection network, which provide advantages for algorithms that make use of a ring structure, is an important issue for multiprocessor systems \cite{Albader}. In this section, we study Hamiltonian cycle embedding of $SSQ_n$ and $BSQ_n$ respectively, implying $SSQ_n$ and $BSQ_n$ are both Hamiltonian.

\begin{theorem}{\bf.}\label{Hamil-SSQn} There exists a Hamiltonian cycle in $SSQ_n$ for $n>2$.
\end{theorem}

\noindent{\bf Proof.} We construct a Hamiltonian cycle $H_1$ of $SSQ_n$ for $n>2$. Start $H_1$ with the vertex $u=\underbrace{0000}\underbrace{0000}\cdots\underbrace{0000}00$ by changing $u_4^1$. Since each vertex is contained in exactly one $K_4$ of $SSQ_n$, it is easy to find a Hamiltonian path $P_1$ of this $K_4$ (vertices with the first 4-bit 00**) from $u$ to $u_1=\underbrace{0011}\underbrace{0000}\cdots\underbrace{0000}00$, where ``*'' means the bit we do not care. Note that $u_1'=\underbrace{1100}\underbrace{0000}\cdots\underbrace{0000}00$ is a neighbor of $u_1$, and $u_1$ is contained in exactly one $K_4$ with the first 4-bit 11**. Similarly, there exists a Hamiltonian path $P_2$ of this $K_4$ (vertices with the first 4-bit 11**) from $u_1'$ to $u_1''=\underbrace{1111}\underbrace{0000}\cdots\underbrace{0000}00$. Concatenating $P_1$ and $P_2$ yields a path $P$ with eight vertices.

We then iteratively change $u_4^j$ ($1\leq j\leq\frac{n-2}{4}-1$) one step as the first 4-bit of the vertices varying on $P$, to reach a new neighbor. In particular, if $j=0$, then we change $u_4^0=u_1u_0$ one step in a cyclic manner with respect to 00, 01, 10, 11. Since we have arrived at a vertex with distinct $u_4^j$, we change its first 4-bit as in the manner above. Consequently, we can obtain a Hamiltonian cycle of $SSQ_n$ for $n\geq2$.

This completes the proof.\qed

To illustrate the approach used in the proof of Theorem \ref{Hamil-SSQn}, we give the following example.

\noindent{\bf{\normalsize Example 1.}} To obtain a Hamiltonian cycle $H_1$ of $SSQ_6$, we starts $H_1$ at $u=000000$ as follows:

{\small
\begin{equation*}
\begin{array}{lll}
H_1: &000000{\rightarrow} 000100{\rightarrow} 001000 {\rightarrow}001100{\rightarrow}110000{\rightarrow}110100{\rightarrow}111000{\rightarrow}111100{\rightarrow}\\
&111101{\rightarrow} 111001{\rightarrow} 110101 {\rightarrow}110001{\rightarrow}001101{\rightarrow}001001{\rightarrow}000101{\rightarrow}000001{\rightarrow}\\
&000010{\rightarrow}000110{\rightarrow}001010{\rightarrow}001110{\rightarrow}110010{\rightarrow}
110110{\rightarrow}111010{\rightarrow}111110{\rightarrow}\\ &111111{\rightarrow}111011{\rightarrow}110111{\rightarrow}110011{\rightarrow}001111{\rightarrow}
001011{\rightarrow}000111{\rightarrow}000011{\rightarrow}\\
&000000
\end{array}
\end{equation*}
}

\begin{lemma}{\bf.}\label{Hamil-BSQ6} There exists a Hamiltonian cycle of $BSQ_6$.
\end{lemma}
\noindent{\bf Proof.} We construct a Hamiltonian cycle $H_2$ of $BSQ_6$ by greedily changing the first two bits of vertices in a cyclic manner with respect to 00, 01, 10, 11 as follows:
{\small
\begin{equation*}
\begin{array}{lll}
H_2: &000000{\rightarrow}010000{\rightarrow}100000{\rightarrow}110000{\rightarrow} 001100{\rightarrow}011100{\rightarrow}101100{\rightarrow}111100{\rightarrow}\\

&001000{\rightarrow}011000{\rightarrow}101000{\rightarrow}111000{\rightarrow} 000100{\rightarrow}010100{\rightarrow}100100{\rightarrow}110100{\rightarrow}\\

&110101{\rightarrow}000101{\rightarrow}010101{\rightarrow}100101{\rightarrow} 111001{\rightarrow}001001{\rightarrow}011001{\rightarrow}101001{\rightarrow}\\

&111101{\rightarrow}001101{\rightarrow}011101{\rightarrow}101101{\rightarrow} 110001{\rightarrow}000001{\rightarrow}010001{\rightarrow}100001{\rightarrow}\\

&100010{\rightarrow}110010{\rightarrow}000010{\rightarrow}010010{\rightarrow}
101110{\rightarrow}111110{\rightarrow}001110 {\rightarrow}011110{\rightarrow}\\

&101010{\rightarrow}111010{\rightarrow}001010{\rightarrow}011010{\rightarrow}
100110{\rightarrow}110110{\rightarrow}000110 {\rightarrow}010110{\rightarrow}\\

&010111{\rightarrow}100111{\rightarrow}110111{\rightarrow}000111{\rightarrow}
011011{\rightarrow}101011{\rightarrow}111011{\rightarrow}001011{\rightarrow}\\

&011111{\rightarrow}101111{\rightarrow}111111{\rightarrow}001111{\rightarrow}
010011{\rightarrow}100011{\rightarrow}110011{\rightarrow}000011{\rightarrow}\\

&000000
\end{array}
\end{equation*}
}\qed

\vskip 0.1 in

By the proof of Lemma \ref{Hamil-BSQ6}, it is known that each path $uvwx$ is contained in $H_2$, where $u_5u_4=00, v_5v_4=01, w_5w_4=10, x_5x_4=11$, and $u_3u_2u_1u_0=v_3v_2v_1v_0=w_3w_2w_1w_0=x_3x_2x_1x_0$. By Remark 2, $u$ and $w$ have the same neighborhood, meanwhile, $v$ and $x$ have the same neighborhood. Since $BSQ_n$ is vertex-transitive, each edge of the cycle $uvwxu$ can be in a Hamiltonian cycle of $BSQ_6$ by modifying $H_2$ through equivalent vertices.

\begin{theorem}{\bf.}\label{Hamil-BSQn} There exists a Hamiltonian cycle in $BSQ_n$ for $n>2$.
\end{theorem}

\noindent{\bf Proof.} For convenience, let $k=\frac{n-2}{4}$. We prove this theorem by induction on $k$. Obviously, by Lemma \ref{Hamil-BSQ6} the statement is true for $k=1$ ($n=6$). Assume that the statement is true for $k-1$ ($k\geq2$). Next we consider $BSQ_n$ ($n=4k+2$). For some $i$ ($1\leq i\leq\frac{n-2}{4}$), all the vertices of the same $u_4^i$ form a subcube $BSQ_{n-4}$. Suppose w.l.o.g. that $i=\frac{n-2}{4}$ and $j=\frac{n-2}{4}-1$. By Definition \ref{variation-2}, there exist sixteen subcubes $BSQ_{n-4}$ with different $u_4^i$s respectively.

\begin{figure}
\centering
\includegraphics[height=80mm]{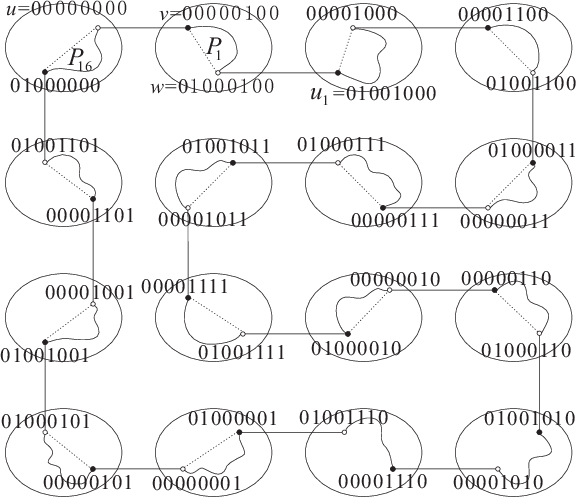}
\caption{Construction of Hamiltonian cycle of $BSQ_n$.} \label{HC-BSQn}
\end{figure}

We construct a Hamiltonian cycle of $BSQ_n$ from $u=000000000000\cdots00$. Since we only consider the first two 4-bits of the vertices in the remainder of the proof, for simplicity, we write $u=00000000$. Clearly, $u$ has a neighbor $v=00000100$ in another $BSQ_{n-4}$. For convenience, let $w=01000100$. By induction hypothesis, this $BSQ_{n-4}$ has a Hamiltonian cycle $H_1$ containing the edge $vw$. Similarly, $w$ has a neighbor $u_1=01001000$. By repeatedly using the above method, a Hamiltonian cycle $H$ of $BSQ_n$ follows (see Fig. \ref{HC-BSQn}), completing the proof.\qed

\vskip 0.1 in

\section{Conclusions}
The shuffle-cube $SQ_n$ was proposed by Li et al. \cite{Li} to enhance some properties of the hypercube. Symmetry is a desirable property in interconnection networks. Highly symmetric network topologies are popular since it often simplifies the computation and routing algorithms. In this paper, we study the symmetric properties of the shuffle-cube $SQ_n$, showing that $SQ_n$ is neither vertex-transitive nor edge-transitive for all $n>2$. To address this shortcoming of the shuffle-cube, we propose two variants of the shuffle-cube, namely simplified shuffle-cube $SSQ_n$ and balanced shuffle-cube $BSQ_n$ and show that they are both vertex-transitive for all $n>2$.

We also design routing schemes for $SSQ_n$ and $BSQ_n$, respectively. Meanwhile, we show that the diameter of $SSQ_n$ and $BSQ_n$ are $\frac{n-2}{2}+2$ and $n$ for $n\geq2$, respectively, which are not greater than that of the hypercube. In addition, we show that both $SSQ_n$ and $BSQ_n$ admit a Hamiltonian cycle, which is an appealing property in interconnection networks.

For a future research, it is meaningful to further determine the automorphism group of the shuffle-cube, and characterize if the simplified shuffle-cube and balanced shuffle-cube are Cayley graphs. In addition, other combinatorial properties of the simplified shuffle-cube and balanced shuffle-cube need to be further investigated, such as fault-tolerant Hamiltonian cycle embedding, edge-disjoint Hamiltonian cycle embedding, super (edge) connectivity, etc.

\vskip 0.3 in

\noindent{\bf\large Appendix}
\vskip 0.2 in
As a by-product of proving vertex-transitivity of $BSQ_n$ in Theorem \ref{Cayley-BSQn}, we mend a flaw in the Property 3 in [IEEE Trans. Comput. 46 (1997) 484--490].
%
%

Property 3 of \cite{Wu} states that for any pair of nodes $v$ and $u$ in $BH_n$ there is an automorphism $T$ of $BH_n$ such that $T(v)=u$, showing that $BH_n$ is vertex-transitive.

Moreover, they showed that if $v=(a_0, a_1,\cdots, a_{n-1})$ and $u=(b_0,b_1,\cdots, b_{n-1})$, then the automorphism $T$ can be defined as:
$$T(w)=(b_0+(-1)^{a_0}(c_0-a_0),b_1+(-1)^{a_0}(c_1-a_1),\cdots,b_{n-1}+(-1)^{a_0}(c_{n-1}-a_{n-1}))$$
for any node $w=(c_0, c_1,\cdots, c_{n-1})$.

In fact, this is not always the case. For example, let $v=(3,0)$ and $u=(1,1)$ in $BH_2$. Then $T(w)=(1+(-1)^3(c_0-3),1+(-1)^3(c_1-0))=(4-c_0,1-c_1)$. We may assume that $x=(0,0)$ and $y=(1,1)$. Then $T(x)=(0,1)$ and $T(y)=(3,0)$. Observe that $xy\in E(BH_2)$, yet $T(x)T(y)\not\in E(BH_2)$ by the definition of $BH_n$. Therefore, the automorphism $T$ of $BH_n$ defined is inappropriate.

In the following, we shall define an automorphism $T$ of $BH_n$ in detail such that $T(v)=u$ for any pair of nodes $v$ and $u$ in $BH_n$.

If $a_{0}$ and $b_{0}$ have the same parity, then
$$T(w)=(-b_{0}+a_{0}+c_{0}, -b_{1}+a_{1}+c_{1},\cdots,-b_{n-1}+a_{n-1}+c_{n-1});$$
If $b_{0}$ is even and $a_{0}$ is odd, then
$$T(w)=(-b_{0}+a_{0}-c_{0}, -b_{1}+a_{1}-c_{1},\cdots,-b_{n-1}+a_{n-1}-c_{n-1});$$
If $b_{0}$ is odd and $a_{0}$ is even, then
$$T(w)=(b_{0}-a_{0}-c_{0}, b_{1}-a_{1}-c_{1},\cdots, b_{n-1}-a_{n-1}-c_{n-1});$$
where addition and substraction are under modulo 4.

Now we are ready to prove that $T$ is an automorphism of $BH_n$.

\begin{theorem}{\bf .}
The function $T$ we have just defined is an automorphism of $BH_n$.
\end{theorem}
\noindent{\bf Proof.} Let $xy$ be an arbitrary edge of $BH_n$ with $x=(x_0, x_1,\cdots, x_{n-1})$ and $y=(y_0, y_1,\cdots, y_{n-1})$. It suffices to show that $T(x)T(y)$ is also an edge of $BH_n$.

If $x$ and $y$ differ in only the first coordinate, then it is straightforward that $T(x)T(y)\in E(BH_n)$. So we assume that $x$ and $y$ differ not only the first coordinate, that is, for some $i\in\{1,\cdots,n-1\}$, $y_i\equiv x_i+(-1)^{x_0}(\text{mod}\ 4)$, and for each $j\in\{1,\cdots,n-1\}\setminus\{i\}$, $y_j=x_j$. The following cases arise.

\noindent{\bf Case 1.} $b_{0}$ and $a_{0}$ are of the same parity. Clearly, we have
$T(x)=(-b_{0}+a_{0}+x_{0}, -b_{1}+a_{1}+x_{1},\cdots,-b_{n-1}+a_{n-1}+x_{n-1})$
and $T(y)=(-b_{0}+a_{0}+y_{0}, -b_{1}+a_{1}+y_{1},\cdots,-b_{n-1}+a_{n-1}+y_{n-1})$.

Thus, the parity of $-b_{0}+a_{0}+x_{0}$ (resp. $-b_{0}+a_{0}+y_{0}$) is the same as that of $x_{0}$ (resp. $y_{0}$). It follows from the definition of $BH_n$ that $T(x)T(y)\in E(BH_n)$.

\noindent{\bf Case 2.} $b_{0}$ is even and $a_{0}$ is odd. Clearly, we have
$T(x)=(-b_{0}+a_{0}-x_{0}, -b_{1}+a_{1}-x_{1},\cdots,-b_{n-1}+a_{n-1}-x_{n-1})$ and
$T(y)=(-b_{0}+a_{0}-y_{0}, -b_{1}+a_{1}-y_{1},\cdots,-b_{n-1}+a_{n-1}-y_{n-1})$. Note that $xy\in E(BH_n)$.

If $x_{0}$ is even, then $y_{i}\equiv x_{i}+1(\text{mod}\ 4)$, implying that $-y_{i}\equiv -x_{i}-1(\text{mod}\ 4)$. Observe that $-b_{0}+a_{0}-x_{0}$ is odd and $-b_{0}+a_{0}-y_{0}$ is even. Combining with $y_j=x_j$ for each $j\in\{1,\cdots,n-1\}\setminus\{i\}$, we have $T(x)T(y)\in E(BH_n)$.

If $x_{0}$ is odd, then $y_{i}\equiv x_{i}-1(\text{mod}\ 4)$, implying that $-y_{i}\equiv -x_{i}+1(\text{mod}\ 4)$. Observe that $-b_{0}+a_{0}-x_{0}$ is even and $-b_{0}+a_{0}-y_{0}$ is odd. Combining with $y_j=x_j$ for each $j\in\{1,\cdots,n-1\}\setminus\{i\}$, we have $T(x)T(y)\in E(BH_n)$.

\noindent{\bf Case 3.} $b_{0}$ is odd and $a_{0}$ is even. In this case, we have
$T(x)=(b_{0}-a_{0}-x_{0}, b_{1}-a_{1}-x_{1},\cdots, b_{n-1}-a_{n-1}-x_{n-1})$
and $T(y)=(b_{0}-a_{0}-y_{0}, b_{1}-a_{1}-y_{1},\cdots, b_{n-1}-a_{n-1}-y_{n-1})$. Observe that $xy\in E(BH_n)$.

If $x_{0}$ is even, then $y_{i}\equiv x_{i}+1(\text{mod}\ 4)$. Thus, $-y_{i}\equiv -x_{i}-1(\text{mod}\ 4)$. Now we see that $b_{0}-a_{0}-x_{0}$ is odd and $b_{0}-a_{0}-y_{0}$ is even. Combining with $y_j=x_j$ for each $j\in\{1,\cdots,n-1\}\setminus\{i\}$, it follows that $T(x)T(y)\in E(BH_n)$.

If $x_{0}$ is odd, then $y_{i}\equiv x_{i}-1(\text{mod}\ 4)$. Thus, $-y_{i}\equiv -x_{i}+1(\text{mod}\ 4)$. Now we see that $b_{0}-a_{0}-x_{0}$ is even and $b_{0}-a_{0}-y_{0}$ is odd. Combining with $y_j=x_j$ for each $j\in\{1,\cdots,n-1\}\setminus\{i\}$, it follows that $T(x)T(y)\in E(BH_n)$.

%
%
%
%

This completes the proof.\qed

\begin{example}{\bf .}
Let $v=(2,0,1,2,3,1)$ and $u=(1,2,3,1,0,2)$ be two vertices in $BH_6$. For any vertex $w$ in $BH_6$, we have
\begin{equation*}
\begin{aligned}
T(w)&=(1-2-c_0,2-0-c_1,3-1-c_2,1-2-c_3,0-3-c_4,2-1-c_5)\\
&=(3-c_0,2-c_1,2-c_2,3-c_3,1-c_4,1-c_5).
\end{aligned}
\end{equation*}
Suppose w.l.o.g. that $x=(x_0,x_1,x_2,x_3,x_4,x_5)$, $y=(x_0\pm1,x_1,x_2,x_3,x_4,x_5)$ and $y'=(x_0\pm1,x_1,x_2+(-1)^{x_0},x_3,x_4,x_5)$. Thus, $xy, xy'\in E(BH_n)$. Moreover, $$T(x)=(3-x_0,2-x_1,2-x_2,3-x_3,1-x_4,1-x_5),$$
$$T(y)=(3-x_0\mp1,2-x_1,2-x_2,3-x_3,1-x_4,1-x_5)$$ and $$T(y')=(3-x_0\mp1,2-x_1,2-x_2-(-1)^{x_0},3-x_3,1-x_4,1-x_5).$$
Clearly, $T(x)T(y), T(x)T(y')\in E(BH_6)$.
\end{example}\qed

\begin{example}{\bf .}
Let $v=(1,2,3,1,2,0,0)$ and $u=(3,1,3,2,1,2,1)$ be two vertices in $BH_7$. For any vertex $w$ in $BH_7$, we have
\begin{equation*}
\begin{aligned}
T(w)&=(1-3+c_0,2-1+c_1,3-3+c_2,1-2+c_3,2-1+c_4,0-2+c_5,0-1+c_6)\\
&=(-2+c_0,1+c_1,c_2,-1+c_3,1+c_4,-2+c_5,-1+c_6).
\end{aligned}
\end{equation*}
Suppose w.l.o.g. that $x=(x_0,x_1,x_2,x_3,x_4,x_5,x_6)$, $y=(x_0\pm1,x_1,x_2,x_3,x_4,x_5,x_6)$ and $y'=(x_0\pm1,x_1,x_2,x_3+(-1)^{x_0},x_4,x_5,x_5)$. Thus, $xy, xy'\in E(BH_n)$. Moreover, $$T(x)=(-2+x_0,1+x_1,x_2,-1+x_3,1+x_4,-2+x_5,-1+x_6),$$
$$T(y)=(-2+x_0\pm1,1+x_1,x_2,-1+x_3,1+x_4,-2+x_5,-1+x_6)$$ and $$T(y')=(-2+x_0\pm1,1+x_1,x_2,-1+x_3+(-1)^{x_0},1+x_4,-2+x_5,-1+x_6).$$
Clearly, $T(x)T(y), T(x)T(y')\in E(BH_7)$.
\end{example}\qed

\end{document}